\documentclass[12pt,bezier]{amsart}
\usepackage{amsfonts,latexsym,amscd}
\newtheorem{theorem}{Theorem}[section]

\newtheorem{proposition}[theorem]{Proposition}

\newtheorem{lemma}[theorem]{Lemma}

\newcommand{\thm}[1]{Theorem \ref{#1}}
\newcommand{\prop}[1]{Proposition \ref{#1}}
\newcommand{\lma}[1]{Lemma \ref{#1}}

\newcommand{\rarrow}{\rightarrow}
\newcommand{\fd}{\operatorname{fd}}

\newcommand{\tor}{\operatorname{Tor}}
\newcommand{\Hom}{\operatorname{Hom}}
\newcommand{\chr}{\operatorname{char}}

\newcommand{\F}{{\mathbb F}}
\newcommand{\calL}{{\mathcal L}}

\newcommand{\calA}{{\mathcal A}}
\newcommand{\Ap}{{\mathcal A}_{\ell}}

\newcommand{\calF}{{\mathcal F}}

\newcommand{\calS}{{\mathcal S}}

\newcommand{\calV}{{\mathcal V}}

\begin{document}

\title[Simplicial Algebras with Vanishing Andr\'e-Quillen Homology]
{On Simplicial Commutative Algebras with Vanishing Andr\'e-Quillen 
Homology}
\author{James M. Turner}
\address{Department of Mathematics\\
Calvin College\\
3201 Burton Street, S.E.\\
Grand Rapids, MI 49546}
\email{jturner@calvin.edu}
\thanks{Research was partially supported by an NSF-NATO postdoctoral 
fellowship and by an NSF grant}
\date{\today}



\begin{abstract}
In this paper, we study the Andr\'e-Quillen homology of simplicial 
commutative $\ell$-algebras, $\ell$ a field, having certain vanishing 
properties. When $\ell$ has non-zero characteristic, we obtain an 
algebraic version of a theorem of J.-P. Serre and Y. Umeda that 
characterizes such simplicial algebras having bounded homotopy groups. 
We further discuss how this theorem fails in the rational case and,
as an application, indicate how the algebraic Serre theorem can 
be used to resolve a conjecture of D. Quillen for algebras of finite 
type over Noetherian rings, which have non-zero characteristic.


\end{abstract}

\maketitle

\section*{Overview}

\subsection*{Algebraic Serre Theorem}

The following topological theorem is due to J.-P. 
Serre \cite{Serre} at the prime 2 and to Y. Umeda \cite{Ume} at odd 
primes.

\bigskip

\noindent {\bf Serre's Theorem.}
{\em Let $X$ be a nilpotent space such that $H_{s}(X;\mathbb{F}_{p}) =0$ for 
$s\gg 0$ and each $H_{s}(X;\mathbb{F}_{p})$ is finite dimensional.
Then the following are equivalent
\begin{itemize} 
\item[1.] $\pi_{s}(X)\otimes \mathbb{Z}/p = 0, \quad s \gg 0$;
\item[2.] $\pi_{s}(X)\otimes \mathbb{Z}/p = 0, \quad s\geq 2$.
\end{itemize}}

\bigskip

In \cite{And,Qui2,Qui3}, M. Andr\'e and D. Quillen constructed the
notion of a homology $D_{*}(A|R;M)$ for a homomorphism $R\to A$ of 
simplicial commutative 
rings, with coefficients in a simplicial $A$-module $M$. These homology groups
can be defined as $\pi_{*}(\calL (A|R)\otimes_{A} M)$ where the 
simplicial $A$-module $\calL (A|R)$ is called the {\bf cotangent 
complex} of $A$ over $R$.

We now propose an algebraic 
analogue of Serre's Theorem for simplicial augmented
$\ell$-algebras. 
To accomplish this 
we will take simplicial homotopy $\pi_{*}(-)$ to be the analogue of 
$H_{*}(-;\mathbb{F}_{p})$ and $H^{Q}_{*}(-)=D_{*}(-|\ell;\ell)$ to be the 
analogue of $\pi_{*}(-)\otimes \mathbb{Z}/p$. 
\bigskip

\noindent {\bf Algebraic Serre Theorem.}
{\em Let $A$ be a homotopy connected (i.e. $\pi_{0}A = \ell$) simplicial 
supplemented commutative $\ell$-algebra, with $\chr{\ell}\neq 0$, such that 
$\pi_{*}A$ is a finite graded $\ell$-module. Then the following are equivalent
\begin{itemize}
   \item[1.] $H^{Q}_{s}(A) = 0, \quad s \gg 0$;
   \item[2.] $H^{Q}_{s}(A) = 0, \quad s\geq 2$.
\end{itemize}}

\bigskip

We shall prove this theorem by following Serre's original 
approach in \cite{Serre}. This will require pooling technical tools 
such  as an analogue of the 
notion of connected covers of spaces and various ways for  
making computations of the homotopy and homology of simplicial 
commutative algebras. 

The algebraic Serre theorem cannot hold in general when the ground 
field has characteristic zero. At the end of \S 2, we indicate a 
partial result in the rational case and point to some examples that 
show that a full version of our theorem cannot hold rationally.

\subsection*{Connections to Quillen's Conjecture}


D. Quillen has conjectured that the cotangent complex has certain 
rigidity properties. In particular, we recall the 
following, which can be found in \cite[(5.7)]{Qui2}:

\bigskip

\noindent {\bf Quillen's Conjecture.} {\em If $A$ is an algebra of 
finite type over a Noetherian ring $R$, such that $A$ has finite flat 
dimension over $R$ and $\fd_{A}\calL(A|R)$ is finite, then $A$ is a 
quotient of a polynomial ring by an ideal generated by a regular 
sequence.}

\bigskip

Earlier results of Lichtenbaum-Schlessinger \cite{LS}, Quillen 
\cite{Qui2}, and Andr\'e \cite{And} prove that an $R$-algebra $A$ is 
a complete intersection if and only if $\fd_{A}\calL (A|R)\leq 
1$. In characteristic 0 the conjecture was proved by 
Avramov-Halperin \cite{AH}. The general case was proved by L. Avramov. 
Furthermore, Avramov characterized 
those homomorphisms $R\to A$ of Noetherian rings having locally finite 
flat dimension with $\fd_{A}\calL(A|R)<\infty$. See \cite{Avr2} for 
details.

As a consequence of the Algebraic Serre Theorem, we have the 
following:

\begin{theorem}\label{quiconj}
Quillen's conjecture holds provided the algebra $A$ has non-zero characteristic.
\end{theorem}

{\it Proof.} Since $\fd_{A}\calL(A|R)\leq N$ if and only if
$D_{s}(A|R;-)=0$ for $s> N$ then 
we seek to show that the latter implies $D_{s}(A|R;-)=0$ for $s\geq 2$. By 
\cite[(S.30)]{And}, it is enough to show that, for each prime ideal 
$\wp\subset A$, $D_{s}(A|R;\mathbf{k}(\wp))=0$ for $s\geq 2$, 
where $\mathbf{k}(\wp)$ is the residue field of $A_{\wp}$. 
Since $A$ has non-zero characteristic then each $\mathbf{k}(\wp)$ has prime 
characteristic. Let $\ell$ denote a fixed residue field.

Since $A$ is an algebra of finite type over $R$, then the unit map 
factors as $R \rightarrow R[X] 
\stackrel{\sigma}{\rightarrow} A$, with $X$ a finite set and $\sigma$ 
a surjection. Since $R\to R[X]$ is a flat homomorphism, then 
$D_{s}(R[X]|R;\ell) \cong D_{s}(\ell[X]|\ell;\ell) = 0$ for $s\geq 1$,
by \cite[(4.54, 6.26)]{And}. 
An application of \cite[(5.1)]{And} now implies that
$D_{s}(A|R;\ell) \cong D_{s}(A|R[X];\ell)$ for $s\geq 2$.  
Since $\fd_{R}A = \fd_{R[X]}A$, by a change-of-rings spectral sequence 
argument, we may thus 
assume that $R\to A$ is surjective.

Let $\calF$ be the 
homotopy pushout over $\ell$ of $R\to A$ in the simplicial model 
category of simplicial commutative $R$-algebras over $\ell$ (see 
\cite{Qui1, Qui2, Qui3, Goe2} for general discussions pertaining to this 
model structure). Then $\calF$ is a connected simplicial supplemented
commutative $\ell$-algebra with the properties
$$
	D_{*}(\calF|\ell;\ell) \cong D_{*}(A|R;\ell)
$$

\noindent and

$$
	\pi_{*}\calF \cong \tor^{R}_{*}(A,\ell),
$$

\bigskip 

\noindent the first isomorphism following from the flat base change property for 
Andr\'e-Quillen homology \cite[(4.7)]{Qui3} while the second follows 
from an argument utilizing the 
Kunneth spectral sequence of Theorem 6.b in \cite[\S II.6]{Qui1}. 

By the assumption that $\fd_{R}A< \infty$, it follows that 
$\pi_{*}\calF$ is a finite graded $\ell$-module. The result 
now follows. \hfill $\Box$

\subsection*{Generalizing Quillen's Conjecture}

We propose the following 
simplicial generalization of Quillen's conjecture. 

\bigskip


\noindent {\bf Conjecture.} {\em Let $R$ be a Noetherian ring and let 
$A$ be a simplicial commutative $R$-algebra with the following properties:
\begin{itemize}
	
	\item[(1)] $\pi_{0}A$ is a Noetherian ring having non-zero 
	characteristic;
	
	\item[(2)] $\pi_{*}A$ is finite graded as a $\pi_{0}A$-module;
	
	\item[(3)] $\fd_{R}\pi_{*}A< \infty$.
	
\end{itemize}
Then $D_{s}(A|R;-)=0$ for $s\gg 0$ implies $D_{s}(A|R;-)=0$ for $s\geq 2$.}

\bigskip

\noindent {\em Note.} The condition on the characteristic of $R$ is clearly 
needed, as noted above.

\bigskip

A proof of this conjecture 
can be given when stronger conditions on $\pi_{0}A$ are assumed. See 
\cite{Tur}. For example, 
by the same reduction to the algebraic Serre theorem performed in 
the proof of \thm{quiconj}, the following special case can be proved.

\begin{theorem}
The conjecture holds if property (1) is replaced by the stronger 
property 
\begin{itemize}
	
	\item[$(1^{\prime})$] $\pi_{0}A$ is an algebra of finite type 
	over $R$.
	
\end{itemize}
\end{theorem}

\subsection*{Organization of this paper}
In the first section, we review the needed notions of the model
category structure of simplicial supplemented commutative algebras. 
In particular, we review the construction and some properties of the
homotopy and Andr\'e-Quillen homology for simplicial commutative
algebras. In the next section, we introduce the notion of n-connected 
envelopes for simplicial commutative algebras which dualizes the notion 
of n-connected covers of spaces. We then pause to record 
a crucial splitting result and discuss specific types of simplicial 
commutative algebras which demonstrate the failure of the algebraic 
Serre theorem rationally. We then, in  
the third section, discuss the properties of the Poincar\'e series 
for the homotopy of a simplicial commutative algebra. This leads to 
the last section where we give a proof of the algebraic Serre theorem.

\subsection*{Acknowledgements}

The author would like to thank Haynes Miller and Paul Goerss for  
several conversations relating to this project, Jean Lannes for his 
generous hospitality while the author was staying in France, as well as
for discussing several areas related to this topic, and Lucho Avramov for 
enlightening the author on many aspects of commutative algebra and 
for reading and commenting on several drafts of this paper. The 
author would also like to thank the referee for helping to effectively 
streamline the presentation contained here.

During the time this and related projects were being worked on, the 
author had been a guest visitor at the I.H.E.S., the Ecole 
Polytechnique, and Purdue University. Many thanks to each of these 
institutions for their hospitality and the use of their facilities.

\setcounter{section}{0}

\section{The homotopy and homology of simplicial commutative algebras}

We now review the closed simplicial model category structure for $s
\calA_{\ell}$ the category of simplicial commutative 
$\ell$-algebras augmented over $\ell$.  We will assume the 
reader is familiar with the
general theory of homotopical algebra given in \cite{Qui1}.

We call a map $f: \ A \rarrow B$ in $s \calA_{\ell}$ a
\begin{itemize}
\item[(i)] weak equivalence ($\stackrel{\sim}{\rightarrow}$)
$\Leftarrow\!\!\Rightarrow \pi_* f$ is an isomorphism;
\item[(ii)] fibration ($\rightarrow\!\!\!\!\rightarrow$)
$\Leftarrow\!\!\Rightarrow f$
provided the induced canonical map $A \to B\times_{\pi_{0}B} \pi_{0}A$ is a 
surjection;
\item[(iii)] cofibration($\hookrightarrow$) $\Leftarrow \!\!\Rightarrow f$
is a retract of an almost free map \cite[p. 23]{Goe1}.
\end{itemize}

\begin{theorem}\label{thm2.1} \cite{Qui1, Mil, Goe1}
 With these definitions, $s \calA_{\ell}$ is a closed simplicial
model category.
\end{theorem}

For a description of the simplicial structure, see section II.1 of 
\cite{Qui1}. The details will not be needed for our purposes. 
Given a simplicial vector space $V$, over a field $\ell$, define its 
normalized chain complex $NV$ by
\begin{equation}
N_nV = V_n/(\mbox{Im} \, s_0 + \cdots + \mbox{Im} \, s_n)
\end{equation}
and $\partial: \ N_n V \rarrow N_{n-1}V$ is $\partial = \sum^n_{i=0}
(-1)^id_i$. The homotopy groups $\pi_* V$ of $V$ is defined as
$$
\pi_n V = H_n(NV), \quad n \geq 0.
$$
Thus for $A$ in $s \calA_{\ell}$ we define $\pi_* A$ as above.
The Eilenberg-Zilber theorem (see \cite{Mac}) shows that the algebra
structure on $A$ induces an algebra structure on $\pi_*A$.

If we let $\calV$ be the category of $\ell$-vector spaces, then there
is an adjoint pair
$$
S: \ \calV \Leftarrow\!\!\Rightarrow \calA_{\ell}: \, I,
$$
where $I$ is the augmentation ideal function and $S$ is the symmetric
algebra functor.  For an object $V$ in $\calV$ and $n \geq 0$, let
$K(V,n)$ be the associated Eilenberg-MacLane object in $s \calV$ so that
$$
\pi_s K(V,n) = \left\{\begin{array}{ll}
V & s = n; \\[2mm]
0 & s \neq n.
\end{array}\right.
$$
Let $S(V,n) = S(K(V,n))$, which is an object of $s \calA_{\ell}$ 
called a {\bf sphere algebra}.

Now recall the following standard result which will be useful for
us (see section II.4 of \cite{Qui1}).

\begin{lemma}\label{lma3.4}
If $V$ is a vector space, $A$ a simplicial commutative algebra, 
and [\quad,\quad] denotes morphisms in $Ho(s\calA_{\ell})$, then the map 
$$
[S(V,n),A] \rarrow \Hom_{\calV}(V, I\pi_n A)
$$
is an isomorphism. In particular, $\pi_{n}A = [S(n),A]$, where 
$S(n) = S(\ell,n)$.
\end{lemma}

\noindent Here $\calV$ is the category of vector spaces.

Thus the primary operational structure for
the homotopy groups in $s \calA_{\ell}$ is determined by $\pi_*
S(V_{\bullet})$ for any $V_{\bullet}$ in $s \calV$.  By Dold's theorem 
\cite{Dold} there
is a triple $\calS$ on graded vector spaces so that
\begin{equation}
\pi_* S(V) \cong \calS(\pi_* V)
\end{equation}
encoding this structure. If $\chr{\ell} = 0$, $\calS$ is the free skew
symmetric functor and, for $\chr{\ell} > 0$, $\calS$ is a certain free divided
power algebra  (see, for example, \cite{Bou, Goe1, Nic}).

Recall \cite{Qui3} that given a map of simplicial commutative rings 
$R\to S$, there is a functorially defined simplicial $S$-module 
$\Omega_{S|R}$ called the {\bf Kaehler differentials of $S$ over $R$}. 
Replacing $S$ by a cofibrant simplicial $R$-algebra model $X$ then 
the {\bf cotangent complex of $S$ over $R$} is defined as the 
cofibrant simplicial $S$-module
$$
\calL (S|R) := \Omega_{X|R}\otimes_{X}S
$$
and the {\bf Andr\'e-Quillen homology of $S$ over $R$} with 
coefficients in a simplicial $S$-module $M$ is defined as
$$
D_{*}(S|R;M) := \pi_{*}(\calL (S|R)\otimes_{S}M).
$$

For $A$ in $\calA_{\ell}$, define the indecomposable functor to be $QA =
I(A)/I^2(A)$ which is an object of $\calV$.  Define the 
homology functor $H^{Q}_{*}(-): s\calA_{\ell} \to gr\calV$ 
\cite{Goe1,Goe2,Mil} by
$$
H^Q_s(A) = \pi_s QX, \quad s \geq 0,
$$
where we choose a factorization
$$
\ell \hookrightarrow X \stackrel{\sim\hspace*{5pt}}{\rightarrow \!\!\!\!\!
\rightarrow} A
$$
of the unit $\ell \rarrow A$ as a cofibration and a trivial
fibration.  This definition is independent of the choice of
factorization as any two are homotopic over $A$ (note that every object
of $s \calA_{\ell}$ is fibrant). It is straightforward to show 
\cite[(A.1)]{Goe1} that 
$$
\Omega_{B|\ell}\otimes_{B}\ell \cong QB,
$$
for any augmented $\ell$-algebra $B$, and so
$$
H^Q_{*}(A) = D_{*}(A|\ell;\ell).
$$

We now summarize methods for computing homotopy and Andr\'e-Quillen
homology that we will need for this paper.

\begin{proposition}\label{prop2.4}
\begin{itemize}
\item[(1)] If $f: \ A \stackrel{\sim}{\rarrow} B$ is a weak equivalence
in $s \calA_{\ell}$, then $H^Q_*(f): \ H^Q_*(A)
\stackrel{\cong}{\rarrow} H^Q_*(B)$ is an isomorphism. The converse 
holds provided $I\pi_{0}A = 0$, that is, $A$ is {\bf homotopy 
connected}. 
\item[(2)] There is a Hurewicz homomorphism $h: \ I\pi_* A \rarrow
H^Q_*(A)$ such that if $A$ is homotopy connected and 
$H^Q_s (A) = 0$ for $s<n$ then $A$ is $(n-1)$-connected and
\begin{itemize}
  \item[i.] $h: \ \pi_n A \stackrel{\cong}{\rarrow} H^Q_n(A)$ is an 
    isomorphism and
  \item[ii.] $h: \ \pi_{n+1} A \stackrel{\cong}{\rarrow} H^Q_{n+1}(A)$ 
    is a surjection, which is also injective for $n>1$. 
\end{itemize}
\item[(3)] Let $A \stackrel{f}{\rarrow} B \stackrel{g}{\rarrow} C$ be a 
cofibration sequence in $Ho(s \calA_{\ell})$.  Then:
There is a long exact sequence
$$
\begin{array}{l}
\cdots \rarrow H^Q_{s+1} (C) \stackrel{\partial}{\rarrow} H^Q_s(A)
\stackrel{H^Q_*(f)}{\rarrow} H^Q_s(B) \\[3mm]
\hspace*{20pt} \stackrel{H^Q_*(g)}{\rarrow} H^Q_s(C)
\stackrel{\partial}{\rarrow} H^Q_{s-1}(C) \rarrow \cdots
\end{array}
$$
\end{itemize}
\end{proposition}

{\em Proof.}
For all of these, see \cite[IV]{Goe1}.  In particular, (2) follows 
from Quillen's fundamental spectral sequence and the connectivity of Dold's 
functor $\calS$ \cite{Dold}. 
\hfill $\Box$
\bigskip


\section{Connected Envelopes}

In this section, we construct and determine some
properties of a useful tool for studying simplicial algebras.

Given $A$ in $s\Ap$, which is homotopy connected, we define 
its {\bf connected envelopes} to be a sequence of cofibrations
$$
A = A(0) \stackrel{j_{1}}{\rightarrow} A(1) \stackrel{j_{2}}{\rightarrow} \cdots
\stackrel{j_{n}}{\rightarrow} A(n) \stackrel{j_{n+1}}{\rightarrow} \cdots
$$
with the following properties: 
\begin{itemize}
\item[(1)] For each $n\geq 1$, $A(n)$ is a $n$-connected. 
\item[(2)] For $s > n$,
$$
H^Q_s A(n) \cong H^Q_sA.
$$
\item[(3)] There is a cofibration sequence
$$
S(H^Q_{n} A, n) \stackrel{f_{n}}{\rarrow} A(n-1) \stackrel{j_{n}}{\rarrow} A(n).
$$
\end{itemize}

The existence of a connected envelopes is a consequence of the 
following :

\begin{proposition}\label{prop3.5}
Let $A$ in $s \calA_{\ell}$ be $(n-1)$-connected for $n \geq 1$.
Then there exists a map in $s \calA_{\ell}$,
$$
f_n: \ S(H^Q_nA,n) \rarrow A,
$$
with the following properties
\begin{itemize}
\item[1.] $f_{n}$ is an isomorphism on $\pi_n$ and $H^Q_n$;
\item[2.] the homotopy cofibre $M(f_n)$ of $f_n: \
S(H^Q_nA,n) \rarrow A$ is $n$-connected and satisfies $H^Q_sM(f_n) \cong
H^Q_sA$ for $s > n$; 
\item[3.] if $H^Q_s A = 0, \
s \neq n > 0$ then $f_{n}$ is an isomorphism in $Ho(s \calA_{\ell})$.
\end{itemize}
\end{proposition}

{\em Proof.} 
(1.) By the Hurewicz theorem, \prop{prop2.4} (2), the map
$h: \ \pi_nA \rarrow H^Q_n A$ is an isomorphism. 
By \lma{lma3.4} we have an isomorphism
\begin{eqnarray*}
[S(H^Q_nA,n),A]
& \cong & \Hom_{\calV}(H^Q_n A, I\pi_n A).
\end{eqnarray*}
Choosing $f_n$ to correspond to the inverse of $h$
gives the result.

(2.) This follows from (1.) and the transitivity sequence 
$$
H^Q_{s+1}M(f_n) \rarrow H^Q_s S(H^Q_n A,n) \rarrow H^Q_s A \rarrow
H^Q_s M(f_n).
$$

(3.) By (1.), $f_n: \ S(H^Q_nA,n) \rarrow A$ is
an $H^Q_n$-isomorphism and hence a weak equivalence by
\prop{prop2.4}(1). The converse follows from the computation
$$
H^Q_s S(V,n) = \pi_s QS(V,n) = \pi_s K(V,n) = V
$$
for $s = n$ and 0 otherwise.
\hfill $\Box$
\bigskip

\subsection*{Applications}

\begin{proposition}\label{thm1.4}
If there is a cofibration sequence in $s\calA_{\ell}$
$$
S(V,n-1) \to A \to S(W,n)
$$
for some vector spaces $V$ and $W$ and some $n>1$, then in 
$Ho(s\calA_{\ell})$
$$
A \cong S(H^{Q}_{n-1}A,n-1)\otimes S(H^{Q}_{n}A,n).
$$
\end{proposition}

{\em Proof.}
\prop{prop2.4} (3) tells us that $H^{Q}_{s}A=0$ for $s\neq n, 
n-1$, and there is an exact sequence
$$
0 \to H^{Q}_{n}A \to V \to W \to H^{Q}_{n-1}A \to 0.
$$
Thus $A$ is $n-2$ connected and a connected envelope gives a 
cofibration
$$
S(H^{Q}_{n-1}A,n-1) \stackrel{i}{\rightarrow} A 
\stackrel{j}{\rightarrow} S(H^{Q}_{n}A,n)
$$
for which $H^{Q}(j)$ is an isomorphism. \lma{lma3.4} and 
\prop{prop2.4} (2) give a commutative diagram
$$
\begin{array}{ccc}
[S(H^Q_nA,n),A] & \stackrel{j_*}{\longrightarrow} &
[S(H^Q_nA,n),S(H^Q_nA,n)] \\[1mm]
\cong \downarrow \hspace*{10pt}
&&
\hspace*{10pt} \downarrow \cong \\[1mm]
\Hom(H^Q_n A,\pi_{n}A) & \stackrel{j_{*}}{\longrightarrow} & \Hom(H^Q_n A, 
\pi_{n}S(H^Q_n A,n)) \\[1mm]
\hspace*{10pt} \cong \downarrow h_{*} \hspace*{0pt}
&&
\hspace*{0pt} h_{*} \downarrow \cong \\[1mm]
\Hom(H^Q_n A, H^{Q}_n A) & \stackrel{\cong}{\longrightarrow}
& \Hom(H^Q_n A, H^Q_n A)
\end{array}
$$
which shows that $j$ splits up to homotopy. \hfill $\Box$

\bigskip

From \prop{prop3.5} (3), if $\chr{\ell} = 0$ and $V$ 
finite-dimensional then $H^{Q}_{*}S(V,n) \cong V$ concentrated in 
degree $n$ and $\pi_{*}S(V,n)$ is free
skew-commutative on a basis of $V$ concentrated in degree $n$. Thus
$\pi_{*}S(V,n)$ is bounded for any odd $n$, showing that the 
algebraic Serre theorem cannot hold rationally. On the other hand, we 
do have the following

\begin{proposition}\label{ratcase}
Let $A$ be a connected simplicial augmented commutative 
$\ell$-algebra, with $char \, \ell = 0$, such that $\pi_{*}A$ is a 
finite graded $\ell$-module. Then if $H^{Q}_{odd}=0$ and 
$H^{Q}_{s}A=0$ for $s\gg 0$ we can conclude that $I\pi_{*}A = 0$.
\end{proposition}

\bigskip

{\em Proof.} 1. Suppose $H^{Q}_{m}A\neq 0$ implies that $2r\leq m \leq 2s$. 
Then $H^{Q}_{m}A(2r)\neq 0$ implies that $2(r+1)\leq m \leq 2s$. 
Furthermore, $H^{Q}_{odd}A = 0$ and $\pi_{*}A(2r)$ is a finite 
graded $\ell$-module by a spectral sequence argument \cite[\S II.6]{Qui1}
(using the fact that $\pi_{*}S(V,2r)$ is finitely-generated 
polynomial, when $V$ is finite). The result follows by an induction 
on $s-r$, given that the result is certainly true for $A<r,1>\simeq 
S(2r)$, for any $r$.   \bigskip $\Box$

\bigskip

{\em Example.}
Here is another example of a type of rational simplicial algebra with 
finite homotopy and Andr\'e-Quillen homology.

Since $\pi_{*}S(2r)\cong \ell[x_{2r}]$, let $f : S(2rs) \to S(2r)$ 
represent $x_{2r}^{s}$. Define $A<r,s>$ to be the cofibre of $f$. 
Then the cofibration sequence extends to 
$$
S(2r) \to A<r,s> \to S(2rs+1).
$$
The computation of $\pi_{*}(A<r,s>)$ can be achieved by a Serre spectral 
sequence argument (see the proof of \lma{poiprop})  and the computation of 
$H^{Q}_{*}(A<r,s>)$ can be obtained from \prop{prop3.5} (3) using 
\prop{prop2.4} (3). In the end, we obtain 
$$
\pi_{m}(A<r,s>) =\\[3mm]
 \begin{cases}
   \ell & \quad m= 2ri, \; 0 \leq i < s, \\
           0 & \quad otherwise
 \end{cases}
$$
and
$$
H^{Q}_{m}(A<r,s>) =\\[3mm]
 \begin{cases}
   \ell & \quad m= 2r, \\
   \ell & \quad m=2rs+1, \\
           0 & \quad otherwise.
 \end{cases}
$$

\section{The Poincar\'e Series of a Simplicial Algebra}

Let $A$ be a homotopy connected simplicial supplemented commutative
$\ell$-algebra such that $\pi_{*}A$ is of finite-type. We define its
{\em Poincar\'e series} by
$$
\vartheta(A,t) = \sum_{n\geq 0}(\dim_{\ell}\pi_{n}A)t^{n}.
$$
If $V$ is a finite-dimensional vector space and $n>0$ we write
$$
\vartheta(V,n,t) = \vartheta(S(V,n),t).
$$
Combining the work of \cite{Car} with \cite{Serre, Ume}, 
this latter series converges in the open unit disc.

Given power series $f(t) = \sum a_{i}t^{i}$ and $g(t) = \sum b_{i}t^{i}$
we define the relation $f(t) \leq g(t)$ provided $a_{i}\leq
b_{i}$ for each $i\geq 0$.

\begin{lemma}\label{poiprop}
Given a cofibration sequence
$$
A \rightarrow B \rightarrow C
$$
of connected objects in $\calA_{\ell}$ with finite-type homotopy
groups, then
$$
\vartheta(B,t) \leq \vartheta(A,t)\vartheta(C,t)
$$
which is an equality if the sequence is split.
\end{lemma}

{\em Proof.} First, there is a Serre spectral sequence
$$
E^{2}_{s,t}=\pi_{s}(C\otimes\pi_{t}A) \Longrightarrow \pi_{s+t}B.
$$
This follows from Theorem 6(d) in $\S$II.6 of \cite{Qui1}, which gives
a $1^{st}$-quadrant spectral sequence
$$
E^2_{*,*} = \pi_* (B \otimes_A \pi_* A) \Rightarrow \pi_* B,
$$
where $\pi_*A$ is an $A$-module via the augmentation $A \rarrow \pi_0
A$.  Here we can assume our cofibration sequence is a cofibration with
cofibre $C$.  Since $A$ is connected, then $B \otimes_A \pi_* A \cong C
\otimes \pi_* A$.  

Thus we have
$$
\vartheta(A,t)\vartheta(C,t) =
\sum_{n}(\sum_{i+j=n}\dim_{\ell}E^{2}_{i,j})t^{n} \geq
\vartheta(B,t).
$$
Finally, if the cofibration sequence is split then the spectral sequence
collapses, giving an equality. \hfill $\Box$
\bigskip

Now given two power series $f(t)$ and $g(t)$ we say $f(t) \sim g(t)$
provided $\lim_{t\to \infty }f(t)/g(t)\\ = 1$. Given a Poincar\'e series
$\vartheta(V,n,t)$, for a finite-dimensional $\ell$-vector space $V$
and $n>0$, let
$$
\varphi(V,n,t) = \log_{p}\vartheta(V,n,1-p^{-t}).
$$
Then the following is a consequence of Th\'eor\`eme 9b in
\cite{Serre} and its generalization to arbitrary non-zero 
characteristics in \cite{Ume}, utilizing the results of \cite{Car} to 
translate into our present venue.

\begin{proposition}\label{poieq}
For $V$ an $\ell$-vector space of finite dimension
$q$ and $n>0$ then $\varphi(V,n,t)$ converges on the real line and
$$
\varphi(V,n,t)\sim qt^{n-1}/(n-1)!.
$$
\end{proposition}

\section{Proof of the Algebraic Serre Theorem.}

Recall that $A$ is to be a connected simplicial augmented commutative 
$\ell$-algebra 
with $H^{Q}_{*}(A)$ bounded and $\pi_{*}A$ a finite graded 
$\ell$-module.
The approach we take is to mimic the proof of Serre's Theorem in 
\cite{Serre}; utilizing higher connected envelopes, in place of higher 
connected covers, and Poincar\'e series for homotopy, in place of 
Poincar\'e series for homology. Unfortunately, owing to the nature of 
cofibration sequences, Serre's original proof runs into a glitch at 
the start in our situation. Fortunately, if we skip the first step and 
evoke \prop{thm1.4}, the remainder of Serre's proof works without a 
hitch.

\bigskip 

\noindent {\em Proof of the Algebraic Serre Theorem.}  
Let 
$$
n = \max\{ s| H^{Q}_{s}(A) \neq 0 \}.
$$ 
We must show that $n = 1$.

Consider the connected envelope
$$
S(H^{Q}_{s}(A),s) \rightarrow A(s-1) \rightarrow A(s)
$$
for each s. From the theory of cofibration sequences (see section I.3 of
\cite{Qui1}) the above sequence extends to a cofibration sequence
$$
A(s-1) \rightarrow A(s) \rightarrow S(H^{Q}_{s}(A),s+1).
$$
Thus, by \lma{poiprop}, we have
$$
\vartheta(A(s),t) \leq
\vartheta(A(s-1),t)\vartheta(H^{Q}_{s}(A),s+1,t).
$$
Starting at $s=n-2$ and iterating this relation, we arrive at the
inequality
$$
\vartheta(A(n-2),t) \leq \vartheta(A,t)\prod_{s=1}^{n-2}
\vartheta(H^{Q}_{s}(A),s+1,t).
$$
Now, $A(n-1) \cong S(H^{Q}_{n}(A),n)$ by \prop{prop3.5} (3),
but, by \prop{thm1.4} (1) and \lma{poiprop}, we have
$$
\vartheta(A(n-2),t) =
\vartheta(H^{Q}_{n-1}(A),n-1,t)\vartheta(H^{Q}_{n}(A),n,t).
$$
Since $\pi_{*}(A)$ is of finite-type and
bounded then there exists a $D>p$ such that
$\vartheta(A,t) \leq D$, in the open unit disc. Combining, we have
$$
\vartheta(H^{Q}_{n-1}(A),n-1,t)\vartheta(H^{Q}_{n}(A),n,t)
\leq D\prod_{s=1}^{n-2}\vartheta(H^{Q}_{s}(A),s+1).
$$
Applying a change of variables and $\log_{p}$ to the above inequality, we get
$$
\varphi(H^{Q}_{n-1}(A),n-1,t)+\varphi(H^{Q}_{n}(A),n,t)
\leq d+\sum_{s=1}^{n-2}\varphi(H^{Q}_{s}(A),s+1).
$$
By \prop{poieq}, there is a polynomial $f(t)$ of degree $n-2$, a
non-negative integer $a$, and positive integers $b$ and $d$ such that
$$
at^{n-2}+bt^{n-1} \leq d+f(t), \; t\gg 0
$$
which is clearly false for $n>1$. Thus $n=1$. The rest of the proof
follows from \prop{prop3.5} (3). \hfill $\Box$

\end{document}